# Mixed-Integer Linear Programming (MILP) for Garment Line Balancing

Prof Dr Ray Wai Man Kong[1], Ding Ning[2], Theodore Ho Tin Kong[3]
[1]Adjunct Professor, System Engineering Department, City University of Hong Kong
[1]Chairman of the Intelligent Manufacturing Committee of the Doctors Think Tank Academy
[1]Director, Eagle Nice International (Holding) Ltd.
[1]VP CityU Engineering Doctorate Society
Hong Kong, China
dr.raykong@cityu.edu.hk, dr.raykong@gmail.com
[2]Engineering Doctorate Student, System Engineering Department
[2]City University of Hong Kong
Hong Kong, China
[3]Graduated Student, Master of Science in Aeronautical Engineering
[3]Hong Kong University of Science and Technology
Hong Kong, China

**Abstract**

This applied research article explores the application of Mixed-Integer Linear Programming (MILP) to address line-balancing challenges in the garment industry, focusing on optimizing production processes under multiple constraints. By integrating MILP with Lean Methodology principles, the study demonstrates significant improvements in operational efficiency and cost-effectiveness. The case study, conducted in collaboration with Prof Dr Ray WM Kong, highlights the successful implementation of MILP using IBM CPLEX Studio to optimize production order quantities across online and offline operations. The results reveal a remarkable reduction in labour costs, exceeding 50%, while effectively managing resource capacity and demand constraints. This study not only validates the theoretical underpinnings of MILP in resolving line-balancing issues but also underscores its practical applicability in modernizing garment production. The findings contribute valuable insights into the potential of advanced optimization techniques to enhance competitiveness and sustainability in the garment industry. This abstract succinctly captures the essence of the research, emphasizing the methodology, results, and significance of the study.

## I. INTRODUCTION

Line balancing is an essential practice for garment manufacturing factories seeking to optimize their operations. By ensuring an even distribution of tasks, manufacturers can enhance resource utilization, reduce work in progress (WIP), improve production efficiency, maintain quality standards, and adapt to changing market conditions. Ultimately, effective line balancing not only drives operational excellence but also positions garment manufacturers for sustained success in a competitive landscape. To address this issue, line balancing techniques have proven to be effective in optimizing operations without incurring additional costs. By matching the output from each operation and calculating operator capacity, line balancing ensures more efficient utilization of resources. This paper aims to provide a comprehensive framework for line balancing in sewing assembly lines within the garment industry, focusing on the intelligent sewing hanger method. By offering valuable insights and practical guidelines, the study seeks to empower manufacturers to analyze their systems, enhance efficiency, and ultimately maximize output.

As the garment industry continues to evolve, adopting innovative approaches to line balancing will be crucial for meeting the demands of a rapidly changing market and ensuring long-term success.

## II. GARMENT LINE BALANCING PROBLEM

### A. Problem of Line Unbalancing in Garment Manufacturing

The Garment Line Balancing Problem is a multifaceted challenge that requires careful consideration of various factors, including task characteristics, resource constraints, and production goals. Garment planning problems encompass a wide range of challenges that extend beyond line balancing. Addressing these issues requires a holistic approach that integrates demand forecasting, production scheduling, inventory management, and quality control, among other factors. By effectively managing these planning problems, garment manufacturers can enhance operational efficiency, reduce costs, and improve customer satisfaction, ultimately leading to greater competitiveness in the market.

By employing systematic approaches and leveraging modern technologies, garment manufacturers can effectively address line balancing issues, leading to improved efficiency, reduced lead times, and enhanced overall productivity. As the industry continues to evolve, the importance of effective line balancing will only grow, making it a critical area for ongoing research and development.

A production line is not balanced; hence, there would be the following production problems:

- More Accumulate WIP:
  Some operations can produce more, and some can produce less, which will increase the production line's Work In Progress (WIP).
- Reduced Efficiency:
  In an imbalanced assembly line, the flow of input and output is uneven. It means that an upstream operation output is a downstream operation input. Because of this reason, some worker will not get loading input as per their capacity of producing output, hence they will be underutilized. In this case, it is to make matters worse more machines and manpower will be allocated to increase production, but efficiency will fall even more.
- Chaos on the production floor:
  Front-line management and workers push themselves to produce more work in process as the chaos of non-bottleneck operations at the imbalanced production line with no results because without improving the line balance all the other efforts will be wasted.
- Production Planning Problem: Determining the optimal schedule for producing different garment styles within a given timeframe can be complex, especially when considering setup times, machine availability, and labour constraints.
- Impact: Inefficient scheduling can lead to delays, increased lead times, and reduced responsiveness to market changes.

### III. LITERATURE REVIEW

Prof Dr Ray Wai Man Kong [1] discusses strategies for reducing Standard Applied Minutes (SAM) and balancing the capacities of machines, machine centers, and work centers at the initial stage of line balancing for output rates. He emphasizes that simply increasing the capacity of individual machines and assembly lines does not necessarily enhance overall garment production output and productivity due to issues related to line imbalance. The article "Lean Methodology for Lean Modernization" outlines a methodology for implementing lean technology to develop a future state of value stream mapping (VSM) and establish goals, while also identifying bottlenecks in the garment manufacturing process to improve capacity and achieve a balanced production workflow.

Referring to Ocident Bongomin [2], Assembly Line Balancing Problem (ALBP) also known as assembly line design, is a family of combinatorial optimization problems that have been widely studied in literature due to its simplicity and industrial applicability. ALBP is an NP-hard as it subsumes the bin packing problem as a special case. ALBPs arise whenever an assembly line is configured, redesigned, or adjusted.

Published literature shows that the scope of the ALBP in research is indeed quite clear, with well-defined sets of assumptions, parameters, and objective functions. However, these borders are frequently transgressed in real-life situations, in particular for complex assembly line systems like most garment manufacturing. The applied line-balancing problems in garment manufacturing evolved because garment assembly line poses unique balancing problems to those of large body assembly lines such as trucks, buses, aircraft, and machines.

It consists of distributing the total workload for manufacturing any unit of the products to be assembled among the workstations along the line subject to a strict or average cycle time. The general principles of line balancing are (1) industrial environments for which the line balancing problems considered are machining, assembly, and disassembly; (2) number of product models: single-model lines, mixed-model lines, multimodal lines; (3) line layout: basic straight line, straight lines with multiple workplaces, U-shaped lines, lines with circular transfer.

The assembly line balancing (ALB) problem has been studied by enterprises for many decades by Gary Yu-Hsin Chen [3]. The ALB model ensures that the staff assignment balances the whole production process to effectively reduce production time or idle time. To meet the ALB, employees' mastery of skills at each task would be considered as an indicator.

However, there are few studies investigating multifunctional (multitasking) workers with multiple levels of skills working at workstations. Our research incorporates the concept of the Toyota Sewing System (TSS) derived from the Toyota Production System (TPS) for the clothing or footwear industry. TSS is credited with less floor space, flexibility and a better working environment. TSS is featured with a U-shaped assembly line and teams of workers making garments on a single-piece flow basis.

Chen et al. [4] address a multi-skill project scheduling problem for IT product development. In their research, the project is divided into multiple projects which are completed by a skilled employee. To solve the scheduling problem, they proposed a multi-objective nonlinear Mixed-Integer programming model. Their research takes into consideration employees' skill proficiency at performing tasks, multifunctional employees and cell formation to minimize the production cycle time. Also, adopt another manner to calculate the cycle time different from the previous studies and further consider the workers' skills to reflect the real-world situation. They find that the production time can be effectively reduced with better personnel assignment and a preferred mode of production system. The human factor is an uncertainty to affects the actual cycle time. It is clarified that is the human factor for actual output and driving the real-time dynamic line balancing of garment assembly.

Hoa Nguyen Thi Xuan advised the Applying Genetic Algorithm for Line Balancing Problem in Garment manufacturing and mentioned that Muhammad Babar Ramzan (2019) used a time study approach to balance the line and improve productivity with results in a 36% increase in machine productivity, reduction of work in process and visibility of the processes also improved. Haile Sime & Prabir Jana (2018) [5] used Arena simulation software to prove the use of simulation techniques in designing and evaluating different alternative production systems from which the one with the best performance can be selected for final implementation. This will help apparel industries to optimize the utilization of their resources through effective line balancing. Markus Proster & Lothar Marz (2015) have shown that dynamic balancing is crucial for high productivity in mixed–model assembly lines to handle the different assembly times of the variants. Common possibilities to treat the resulting capacity peaks are drifting and the allocation of jumpers. A simulation tool was shown that can simulate and visualize these methods and therefore reduce complexity and raise transparency in the planning of assembly lines.

Ghosh and Gagnon (1989) as well as Erel and Sarin (1998) provided detailed reviews on these topics. Configurations of assembly lines for single and multiple products could be divided into three production line types, single–model, mixed–model and multi–model. Single–model assembles only one product, and mixed–model assembles multiple products, whereas a multi-model produces a sequence of batches with intermediate setup operations (Becker & Scholl, 2006).

## IV. METHODOLOGY

*A. Problem of Line Unbalancing in Garment Manufacturing*

Referring to Prof Dr Ray Wai Man Kong's articles, Lean Methodology for Garment Modernization and Line Balancing in the Modern Garment Industry, industrial engineering and lean study are required to study the whole garment manufacturing process flow.

Studying line balance in garment manufacturing through the lens of Industrial Engineering and Lean Technology involves a systematic approach to optimizing production processes. This approach focuses on improving efficiency, reducing waste, and enhancing overall productivity. Here's how these disciplines contribute to the study of line balancing in the garment industry:

1. Understanding Line Balancing Concepts

- Definition of Line Balancing:
  Line balancing refers to the process of assigning tasks to workstations in such a way that each workstation has an equal amount of work, thereby minimizing idle time and maximizing throughput.
- Importance in Garment Manufacturing:
  In garment manufacturing, where production involves a series of sequential operations (e.g., cutting, sewing, finishing), effective line balancing is crucial for meeting production targets and ensuring timely delivery.

2. Data Collection and Analysis

- Time Studies:
  Industrial engineers conduct time studies to determine the time required for each task in the production process. This data is essential for calculating Standard Applied Minutes (SAM) and understanding task durations.
- Workload Analysis:
  Analyzing the workload of each workstation helps identify imbalances and bottlenecks in the production line. This analysis can involve collecting data on operator performance, machine efficiency, and task completion rates.

3. Application of Lean Principles

- Value Stream Mapping (VSM):
  Lean technology emphasizes the use of Value Stream Mapping to visualize the flow of materials and information throughout the production process. VSM helps identify areas of waste, such as excess inventory, waiting times, and unnecessary movements.
- Elimination of Waste:
  Lean principles focus on eliminating the seven types of waste (overproduction, waiting, transport, extra processing, inventory, motion, and defects). By addressing these wastes, manufacturers can improve line balance and overall efficiency.

4. Task Assignment and Workstation Design

- Heuristic Methods:
  Industrial engineers use heuristic methods to assign tasks to workstations based on criteria such as task duration, precedence relationships, and operator skill levels. Common methods include the Largest Candidate Rule and Ranked Positional Weight Method.
- Workstation Design:
- Designing workstations ergonomically and efficiently is crucial for maintaining worker productivity and comfort. This includes considering the layout, tools, and equipment needed for each task.

5. Simulation and Modeling

- Simulation Tools:
  Using simulation software, industrial engineers can model the production line to test different configurations and task assignments. This allows for the evaluation of potential improvements without disrupting actual production.
- What-If Analysis:
  Simulation enables manufacturers to conduct what-if analyses to assess the impact of changes in task assignments, machine capacities, or production schedules on overall line balance and output.

6. Continuous Improvement

- Kaizen:
  Lean technology promotes a culture of continuous improvement (Kaizen), where teams regularly assess processes and seek incremental improvements. This approach encourages ongoing evaluation of line balance and productivity.
- Feedback Loops:
  Establishing feedback mechanisms allows operators and managers to identify issues in real-time and make adjustments to maintain optimal line balance.

7. Performance Metrics

- Key Performance Indicators (KPIs):
  Industrial engineers establish KPIs to measure the effectiveness of line-balancing efforts. Common KPIs include cycle time, throughput, WIP levels, and defect rates.
- Benchmarking:
  Comparing performance metrics against industry standards or best practices helps identify areas for improvement and sets realistic goals for line balancing.

## B. General Garment Manufacturing Process

Referring to Line Balancing in the Garment industry, before the line balancing for the sewing process, the structure of the garment is separated into two major manufacturing processes.

The first one is the part sewing which seems that individual parts sewing. Garment. Components are the basic sections of garments including top fronts, top backs, bottom fronts, bottom backs, sleeves, collars/neckline treatments, cuffs/sleeve treatments, plackets, pockets, and waistline treatments. A few processes are involved in the buttoning, ironing and other equipment for elastic sewing on garment parts which is counted on the part assembly or part sewing process. In the garment factory, it is called the sub-assembly process. Template sewing is one of the automated processes in the automation of part assembly.

The second one is the final main assembly which gets the part assembly to combine to the finished garment. After the garment has been finished with all related main assembly processes, the last operation is trimming, ironing, packing to the polybag and then packing to the carton box.

Before the line balancing for the sewing process, the structure of the garment is separated into two major manufacturing processes. The first one is the part sewing which seems that individual parts sewing. Garment. Components are the basic sections of garments including top fronts, top backs, bottom fronts, bottom backs, sleeves, collars/neckline treatments, cuffs/sleeve treatments, plackets, pockets, and waistline treatments.

A few processes are involved in the buttoning, ironing and other equipment for elastic sewing on garment parts which is counted on the part assembly or part sewing process. In the garment factory, it is called the sub-assembly process. Template sewing is one of the automated processes in the automation of part assembly. The second one is the final main assembly which gets the part assembly to combine to the finished garment.

After the garment has been finished with all related main assembly processes, the last operation is trimming, ironing, packing to the polybag and then packing to the carton box.

## C. Conveyor Line of Sewing Process for Line Balance of Garment

In the traditional batch production layout, the sub-assembly process and main assembly process are located on the same production floor in the batch garment sewing line as the production batch layout in below Fig 1. There is a typical ALBP that can be applied to various mathematic methods to optimize the line balancing, but the travel time is counted for the bundle batch for transportation between one workstation (sewing station) to another workstation. Bundle batch assembly is not easy to handle on transportation and not easy to visualize any overstock of work in progress at production floor.

In the conveyor line and layout in the main assembly, there is a modernization method and way to reduce the travel time between workstations and improve visual manufacturing and front-line control as shown the Fig. 2. The conveyor line and conveyor line layouts have the benefit of line balance on the main assembly output and enhanced efficiency. The Hanger conveyor layout is applied to intelligent manufacturing for garments. Because it does not use the progressive bundle concept, this style of layout eliminates the previous Work-in-Progress. Allows all the materials for a specific garment to be transferred as a unit to any workstation's sewing machine.

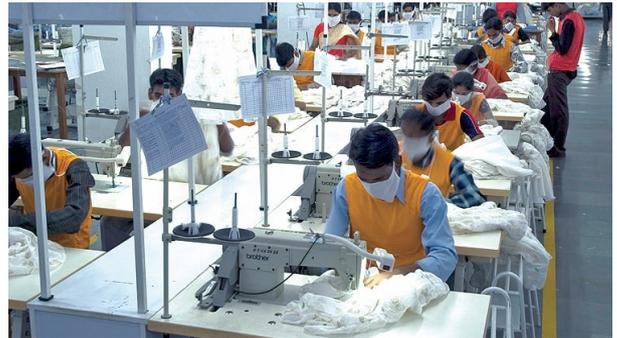

Fig 1 Batch Garment Sewing Line

When an operation at one workstation is completed, the operator should hang the garment to the hanger and then press a button to confirm the finished unit work by faster clipping, so the hanger system can deliver the work-in-progress unit of the garment to the next workstation either mechanically or automatically. It can reduce material handling time. Such a system's layout must be continuous, with no gaps in between. The materials flow through the layout in a loop shape. The hanger line is required to construct the hanger system and equipment. The system is modernized to set up the control device to move the hanger between workstations and provide the just-in-time information to the manufacturing system. The line balancing for the hanger line can be optimized to increase production efficiency by increasing the through-put time based on increased the capacity of the bottleneck workstations in the process as the Lean Methodology for Garment Modernization that Prof Dr Ray WM Kong [6] mentioned with the Design and Experimental Study of Vacuum Suction Grabbing Technology to Grasp Fabric Piece for automation development.

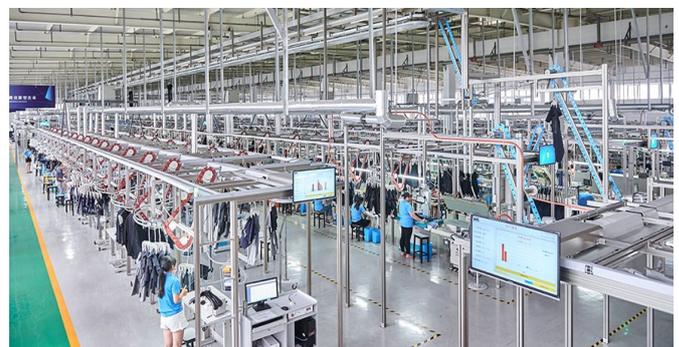

Fig 2 Intelligent Hanger Line and System from INA Intelligent Technology (Zhejiang)

## V. INTEGER PROGRAMMING FOR LINE BALANCING

Garment manufacturing is the most comprehensive process for the line layout design. The prerequisite line balance is required to set up the appropriate line layout for garment manufacture based on various types of garment categories and garment styles. The garment category clarifies the various types of garments: Polo shirts, Dresses, Jeans, Jackets, Pants, Leggings, sportswear, swimwear and others. The garment styles include various garment constructions: grommet drawcord, buttonhole drawcord, inseam gusset, banded hem, banded hem, bound hem exposed trim, elastic of front, zipper of pocket, on-seam pocket and others.

The problem with the fixed facility of the intelligent hanger line is that it cannot be optimized for both part assembly and man flow assembly. The part assembly involves a short cycle time and participating sewing workmanship skills. Optimization and high efficiency are required to reduce the setup time for the part garment assembly in the batch production layout. The operator should continue to produce the same sewing process in the part assembly repeatably. An operator does not change threads, sewing needles and pulling folders if required. The skilled operator can get the benefit of division of work with less change of style, garment construction, fabric piece and upstream sub-assembly work pieces.

*A. Mixed-Integer linear programming (MILP)*

Mixed-Integer programming (MILP) is a mathematical optimization technique that can be effectively used to solve line-balancing problems, particularly in manufacturing environments like garment production. The goal of line balancing is to assign tasks to workstations in such a way that the workload is evenly distributed, minimizing idle time and maximizing throughput. Here's how integer programming can be applied to resolve line-balancing problems:

1. Formulating the Problem

   Define the Decision Variables:

   In integer programming, you first need to define the decision variables. For line balancing, these variables typically represent whether a specific task is assigned to a particular workstation whatever the production quantity.

2. Setting Up the Integer Programming Model

   Once the decision variables, objective function, and constraints are defined, the next step is to set up the integer programming model. This can be done using optimization software or programming languages that support mathematical modelling, such as CPLEX. The study works for the CPLEX for the simulation of the integer programming model.

3. Solving the Integer Programming Model

   After setting up the model, the next step is to solve it using a CPLEX integer programming solver. The CPLEX solver will use algorithms such as branch-and-bound or cutting planes to find the optimal solution that satisfies all constraints while optimizing the objective function subject to the constraint.

4. Interpreting the Results

   Once the solver provides a solution, the results need to be interpreted for the line balancing. Constraints are set up for the integer programming model to optimize the result.

5. Sensitivity Analysis

   After obtaining the optimal solution, it may be beneficial to perform sensitivity analysis to understand how changes in parameters (such as task times or cycle time) affect the solution. For the line balancing, the consumed loading time is the total consumed cycle time with the production quantity. This can help in minimization of the total labour costs (direct cost and indirect cost) in dynamic production environments.

6. Implementation and Continuous Improvement

   Finally, the results from the integer programming model in CPLEX can be implemented in the production environment. Continuous monitoring and feedback can help refine the model and adapt to changes in production requirements, ensuring ongoing efficiency and effectiveness in line balancing.

*B. Mixed-Integer Linear Programming Construction*

The mixed-integer linear programming (MILP) allows for both real and integer variables can be relaxed to a Mixed-Integer linear program by relaxing integer constraints. MILP is often solved using the branch and bound technique. MILP is used in various applications such as task offloading optimization problems.

To minimize labour online cost, $DL_p\,I_p(x)$ and offline cost $IDL_p O_p(x)$ from the direct operator cost in piece rate and indirect labour cost in piece rate representatively, it is the goal the company want to optimize the direct cost and indirect cost. The direct cost is an operator cost to work for the online sewing of garment manufacturing. The indirect cost is related to the supportive staff (line supervisor and additional supportive worker or technician) for working the off-line sewing or joining the online sewing of garment manufacture.

In Table 1, the Notation used in the Mathematical Model defines the various tasks, variable parameters and costs.

Table 1 Notation Used in the Mathematical Models

| Symbol | Definition (unit) |
|---|---|
| Index | |
| *p* | Task Index: $T_1$ .... $T_{19}$ |
| *r* | Resource index of machines (Buttoning MC, Cutting Tools, Manual, Overlock MC, Single Needle MC, Template Sewing MC, Vertical Head Sewing Machine) |
| *n* | Number of Tasks |
| *x* | Number of Production Quantity (pcs) |

| Parameter | |
|---|---|
| **DL** | Direct Labor Cost Rate (RMB/sec) |
| **$DL_p$** | Direct Labor Cost Rate (RMB/sec) in Task $p$ |
| **IDL** | Indirect Labor Cost Rate (RMB/sec) |
| **$IDL_p$** | Indirect Labor Cost Rate (RMB/sec) in Task $p$ |
| **$I_p(x)$** | Online Production Quantity (pcs) in Task $p$ |
| **$O_p(x)$** | Offline Production Quantity (pcs) in Task $p$ |
| **$S_{p,r}$** | Consumed Time (sec) in Task $p$ and Resource $r$ |
| **$R_r$** | Daily Capacity (sec) in Resource $r$ |
| **$D_p$** | Demand quantity (pcs) in Task $p$ |

## Sewing Process Flow for Shirt Manufacturing

| Seq. # | Operation Description |
|---|---|
| | **Collar Making** |
| 1 | Mark Lining |
| 2 | Collar run-stitch |
| 3 | Collar turn & iron |
| 4 | Collar top-stitch |
| 5 | Collar Band Hem |
| 6 | Collar attach to band |
| 7 | Collar trimming, marking & notching |
| 8 | Collar band centre stitch |
| | **Cuff Section** |
| 9 | Cuff hem |
| 10 | Run stitch cuff |
| 11 | Turn cuff |
| 12 | Iron cuffs |
| 13 | Topstitch Cuff |
| | **Pocket Section** |
| 14 | Mark pocket |
| 15 | Pocket mouth iron |
| 16 | Hem pocket |
| 17 | Crease Pocket |
| 18 | Trim pocket |
| | **Front Section** |
| 19 | Mark front for pocket position |
| 20 | Form Button hole placket |
| 21 | Crease B/H placket (single fold) |
| 22 | Topstitch B/H placket |
| 23 | Sew button placket |
| 24 | Attach pocket (1 pocket) |
| 25 | Sew label at placket |

| Seq. # | Operation Description |
|---|---|
| | **Back Section** |
| 26 | Join upper yoke panel |
| 27 | Attach back yoke with back panel |
| 28 | Back yoke topstitch |
| | **Sleeve Section** |
| 29 | Cut sleeve slit at placket position |
| 30 | Notch Sleeves |
| 31 | Iron upper and lower sleeve placket |
| 32 | Attach Plackets |
| 33 | Close lower placket |
| 34 | Close upper plkt & make diamond |
| | **Assembly Section** |
| 35 | Set front & backs & mark neck for collar |
| 36 | Shoulder attach |
| 37 | Shoulder top stitch |
| 38 | Sleeve Attach |
| 39 | Topstitch Armhole |
| 40 | Side Seam |
| 41 | Collar Attach |
| 42 | Collar Close & insert label |
| 43 | Cuff attach & close |
| 44 | Bottom Hem |
| 45 | Button Hole - Front Placket & Collar |
| 46 | Button Attach |

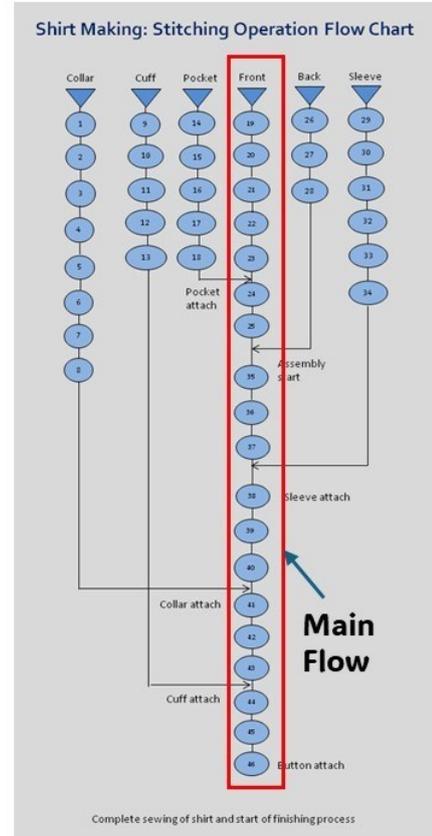

Fig 3 Sewing Process Flow for the Shirt Manufacturing

The minimization formulate is optimized the cost what the company requirement as shown below:

$$Min\ C = \sum_p^n [DL_p I_p(x) + IL_p O_p(x)] \quad (1)$$

The objective of the MILP is to find a solution(s) with a specified quantity with minimum labour costs including direct labour and indirect labour cost. Solutions are considered locally optimized as the principal objective is to find a solution which will define a smooth production by minimizing the objective production quantity balance between online work and offline work.

$$Min\ C = \sum_p^n [DL_p I_p(x) + IL_p O_p(x)]$$

subject to

$$[I_p(x) + O_p(x)] \geq D_p \quad (2)$$
$$\sum_p^n S_{p,r} I_p(x) \leq R_r \quad (3)$$

Where the daily minimization of the sum of all tasks for the total direct labour cost rate multiplied by online production quantity on tasks and the total indirect labour cost rate multiplied by offline production quantity on tasks subject to the sum of online production quantity and offline production quantity by task is greater than the demand quantity in tasks and the sum of consumed loading hour multiplied by online production quantity per task in the resource is greater than the available of the machinery of daily capacity as the supply of resources.

## VI. THE CASE STUDY TO APPLY THE MIXED-INTEGER LINEAR PROGRAMMING FOR LINE BALANCING

The case study to apply the Mixed-Integer linear programming for line balancing in the garment industry proves the theory of Mixed-Integer linear programming to find the minimized total labour costs subject to the daily machinery capacity and demand quantity by tasks.

Referring to the case study of Line Balancing in Modern Garment Industry from Prof Dr Ray WM Kong [1], company X is the shirt manufacturer. The sewing production line selected the sewing process flow for the shirt manufacturing, so the main process includes 19 work tasks as shown the Fig. 3.

In Table 2, the Main Process Sheet of Style A has shown the consumed cycle time to relate the task and 7 resources including machinery and solo manual work.

Based on the demand quantity of 900pcs for style A, the Capacity Requirement Plan of Style A in Table 3 has been calculated for the shortage of online machinery and resources. The offline machinery and resources are planned to compensate for the shortage of machinery and resources.

In the past, the factory manager in company X decided to purchase of shortage of machinery and allocate the offline resources and manpower, but there was no optimization way.

The mixed-integer linear programming and its theory can be combined with the Lean Methodology for the Garment industry and design a new pulling gear for the automation in the International Research Journal of Modernization in Engineering Technology and Science from Prof Dr Ray WM Kong [7], the MILP can calculate the optimized production plan to minimize both online and offline costs of coping with resource and demand constraints. IBM ILOG CPLEX Optimization Studio provides a powerful tool for solving optimization problems, including integer programming (IP) problems. CPLEX uses advanced algorithms and techniques to find optimal or near-optimal solutions to complex mathematical models efficiently. Here's how CPLEX optimizes integer programming problems.

CPLEX primarily uses the branch-and-bound algorithm for solving Mixed-Integer linear programming problems. Initially, CPLEX solves the Mixed-Integer linear programming (MILP) relaxation of the integer programming problem, where the integer constraints are ignored. This provides a bound on the optimal solution in the relaxation stage of the system.

Table 2   Main Process Sheet of Style A

| Task | Description | Resource | CT in resources [after counted workstation (sec/pc)] | | | | | | |
|---|---|---|---|---|---|---|---|---|---|
| | | | Buttoning MC | Cutting Tools | Manual | Overlock MC | Single Sewing MC | Template Sewing MC | Vertical Head Sewing MC |
| T19 | Mark Front for Pocket Position | Manual | 0 | 0 | 30 | 0 | 0 | 0 | 0 |
| T20 | Form Button hole plackets | Cutting Tools | 0 | 40 | 0 | 0 | 0 | 0 | 0 |
| T21 | Crease B/H Placket (Single Fold) | Manual | 0 | 0 | 60 | 0 | 0 | 0 | 0 |
| T22 | Top stitch B/H placket | Single Sewing MC | 0 | 0 | 0 | 0 | 40 | 0 | 0 |
| T23 | Sew Button Placket | Single Sewing MC | 0 | 0 | 0 | 0 | 25 | 0 | 0 |
| T24 | Attach pocket | Manual | 0 | 0 | 20 | 0 | 0 | 0 | 0 |
| T25 | Sewlabel at placket | Template Sewing MC | 0 | 0 | 0 | 0 | 0 | 50 | 0 |
| T35 | Set front & back & mark neck for collar | Manual | 0 | 0 | 60 | 0 | 0 | 0 | 0 |
| T36 | Shoulder attach | Single Sewing MC | 0 | 0 | 0 | 0 | 60 | 0 | 0 |
| T37 | Shoulder top stitch | Single Sewing MC | 0 | 0 | 0 | 0 | 120 | 0 | 0 |
| T38 | Sleeve Attach | Single Sewing MC | 0 | 0 | 0 | 0 | 40 | 0 | 0 |
| T39 | Top stitch armhole | Single Sewing MC | 0 | 0 | 0 | 0 | 80 | 0 | 0 |
| T40 | Side Seam | Overlock MC | 0 | 0 | 0 | 110 | 0 | 0 | 0 |
| T41 | Collar Attach | Single Sewing MC | 0 | 0 | 0 | 0 | 30 | 0 | 0 |
| T42 | Collar Close & Insert Label | Manual | 0 | 0 | 60 | 0 | 0 | 0 | 0 |
| T43 | Cuff Attach & Close | Manual | 0 | 0 | 80 | 0 | 0 | 0 | 0 |
| T44 | Bottom Hem | Vertical Head Sewing MC | 0 | 0 | 0 | 0 | 0 | 0 | 80 |
| T45 | Button Hold - Front Placket & Collar | Manual | 0 | 0 | 70 | 0 | 0 | 0 | 0 |
| T46 | Button Attach (Last Operation for finishing process) | Buttoning MC | 35 | 0 | 0 | 0 | 0 | 0 | 0 |

Table 3   Capacity Requirement Plan of Style A

| Task | Description | Resource | Cycle Time | Loading for 900 orders (sec) | Request no. of M/C (set) | Available no. of M/C (set) | Shortage no. of M/C (Set) |
|---|---|---|---|---|---|---|---|
| T19 | Mark Front for Pocket Position | Manual | 30 | 27,000 | 0.94 | 1 | 0.1 |
| T20 | Form Button hole plackets | Cutting Tools | 40 | 36,000 | 1.25 | 1 | -0.3 |
| T21 | Crease B/H Placket (Single Fold) | Manual | 60 | 54,000 | 1.88 | 2 | 0.1 |
| T22 | Top stitch B/H placket | Single Sewing MC | 40 | 36,000 | 1.25 | 1 | -0.3 |
| T23 | Sew Button Placket | Single Sewing MC | 25 | 22,500 | 0.78 | 1 | 0.2 |
| T24 | Attach pocket | Manual | 20 | 18,000 | 0.63 | 1 | 0.4 |
| T25 | Sewlabel at placket | Template Sewing MC | 50 | 45,000 | 1.56 | 2 | 0.4 |
| T35 | Set front & back & mark neck for collar | Manual | 60 | 54,000 | 1.88 | 2 | 0.1 |
| T36 | Shoulder attach | Single Sewing MC | 60 | 54,000 | 1.88 | 2 | 0.1 |
| T37 | Shoulder top stitch | Single Sewing MC | 120 | 108,000 | 3.75 | 3 | -0.8 |
| T38 | Sleeve Attach | Single Sewing MC | 40 | 36,000 | 1.25 | 1 | -0.3 |
| T39 | Top stitch armhole | Single Sewing MC | 80 | 72,000 | 2.50 | 2 | -0.5 |
| T40 | Side Seam | Overlock MC | 110 | 99,000 | 3.44 | 3 | -0.4 |
| T41 | Collar Attach | Single Sewing MC | 30 | 27,000 | 0.94 | 1 | 0.1 |
| T42 | Collar Close & Insert Label | Manual | 60 | 54,000 | 1.88 | 2 | 0.1 |
| T43 | Cuff Attach & Close | Manual | 80 | 72,000 | 2.50 | 2 | -0.5 |
| T44 | Bottom Hem | Vertical Head Sewing MC | 80 | 72,000 | 2.50 | 2 | -0.5 |
| T45 | Button Hold - Front Placket & Collar | Manual | 70 | 63,000 | 2.19 | 2 | -0.2 |
| T46 | Button Attach (Last Operation for finishing process) | Buttoning MC | 35 | 31,500 | 1.09 | 1 | -0.1 |

** Remark: Red colour means that the shortage of the number of machines and resources

In the CPLEX system, the solution to the MILP relaxation is not integer, CPLEX selects a variable that is fractional and creates two new subproblems (branches) by adding constraints that force the variable to take on integer values (e.g., rounding up or down). It can check the feasibility check whether the current solution satisfies all constraints. If it does, it may update the best-known solution.

CPLEX incorporates various heuristic methods to quickly find feasible solutions, especially for large and complex problems. These heuristics can provide good starting solutions that can be further refined through the branch-and-bound process. Some common heuristics include the algorithms to assign values to variables based on some criteria quickly.

CPLEX finds an optimal or near-optimal solution to cope with the Mixed-Integer linear programming problem, it provides detailed output including the online production quantity, and offline production quantity to minimize the total labour costs.

The CPLEX program source has been created for the Mixed-Integer linear programming for the model formulation (1), (2) and (3) for the line balancing of the shirt sewing process as shown in Fig. 4. The CPLEX software can use the input constraint to calculate the optimized result as minimized both direct and indirect costs.

The result from the CPLEX Studio can debug any programming errors of program source codes. It can show the CPLEX proposed result in Fig 5.

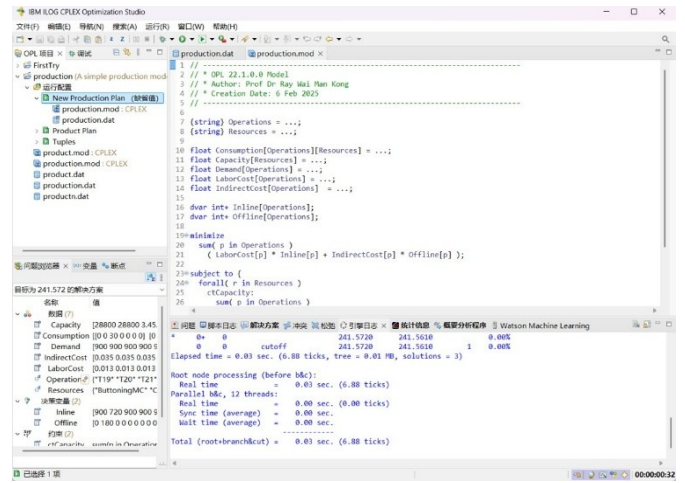

Fig 5 Result from the CPLEX Studio

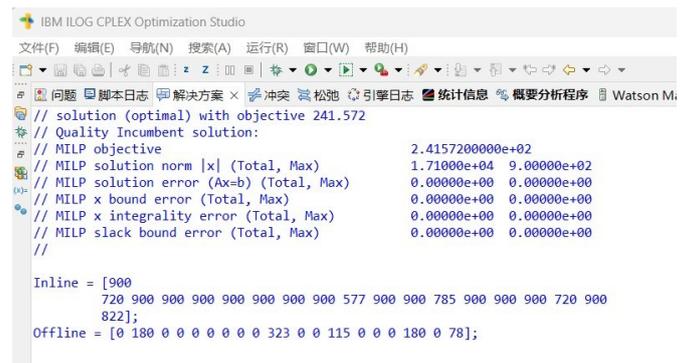

Fig 6 CPLEX Proposed Result

Referring to the CPLEX proposed result in Fig. 6, the MILP objective, online quantity and offline quantity have been calculated as shown in Fig 7 and Fig 8.

```
// ---------------------------------------------------
// * OPL 22.1.0.0 Model
// * Author: Prof Dr Ray Wai Man Kong
// * Creation Date: 6 Feb 2025
// ---------------------------------------------------

{string} Operations = ...;
{string} Resources = ...;

float Consumption[Operations][Resources] = ...;
float Capacity[Resources] = ...;
float Demand[Operations] = ...;
float LaborCost[Operations] = ...;
float IndirectCost[Operations]  = ...;

dvar int+ Inline[Operations];
dvar int+ Offline[Operations];

minimize
  sum( p in Operations )
    ( LaborCost[p] * Inline[p] + IndirectCost[p] * Offline[p] );

subject to {
  forall( r in Resources )
    ctCapacity:
      sum( p in Operations )
        Consumption[p][r] * Inline[p] <= Capacity[r];

  forall(p in Operations)
    ctDemand:
      Inline[p] + Offline[p] >= Demand[p];
}
```

Fig 4 CPLEX programming for the Mixed-Integer Linear Programming

| Operat...大小 19 | 值 |
|---|---|
| T19 | 900 |
| T20 | 720 |
| T21 | 900 |
| T22 | 900 |
| T23 | 900 |
| 24 | 900 |
| T25 | 900 |
| T35 | 900 |
| T36 | 900 |
| T37 | 577 |
| T38 | 900 |
| T39 | 900 |
| T40 | 785 |
| T41 | 900 |
| T42 | 900 |
| T43 | 900 |
| T44 | 720 |
| T45 | 900 |
| T46 | 822 |

Fig 7 Online Production Quantity from CPL

Fig 8 Offline Production Quantity from CPLEX

VII. Conclusion

A. Case Study Analysis

The application of Mixed-Integer Linear Programming (MILP) to resolve line balancing problems in garment production has demonstrated significant potential for optimizing task assignments and enhancing operational efficiency. By systematically formulating the problem with appropriate variables, objective functions, and constraints, manufacturers can effectively align their production processes with customer demand while balancing workloads across both online and offline production environments.

In the case of Style A, the initial production analysis revealed that the hanger line, consisting of 32 seats, was unable to meet the total garment demand of 900 pieces per day, resulting in substantial overstock and understock issues in work-in-progress (WIP) during garment assembly. The throughput rate was recorded at only 720 pieces per day, indicating a clear inefficiency in the production process.

Utilizing Lean Methodology for garment modernization, as outlined by Prof. Dr. Ray Wai Man Kong [8], the original state of the Visual Stream Mapping (VSM) indicated a total daily cost of RMB 589.5. This figure highlighted the financial implications of the existing inefficiencies within the production line.

In contrast, the optimization plan derived from CPLEX's MILP calculations revealed a remarkable reduction in total costs. The optimized total costs were calculated at RMB 241.6 per day, representing a cost saving of over 59% compared to the original state, as demonstrated by the formula: [(RMB 589.5 - RMB 241.6) / RMB 589.5]. This substantial reduction in costs underscores the effectiveness of the MILP approach in minimizing both direct and indirect labour costs associated with garment production.

$$Cost\ Saving\ \% = \frac{(New\ Cost - Old\ Cost)}{Old\ Cost} \times 100\% \quad (4)$$

Table 4 Production Analysis of Style A (Before Line Balancing)

| Task | Description | Resource | Cycle Time (sec/pc) | Number of workstation | Cycle Time after counted workstation (sec/pc) | Daily Output (pcs/day) | Remark |
|---|---|---|---|---|---|---|---|
| T19 | Mark Front for Pocket Position | Manual | 30 | 1 | 30.0 | 960 | Overstock of WIP |
| T20 | Form Button hole plackets | Cutting Tools | 40 | 1 | 40.0 | 720 | Less Stock of WIP |
| T21 | Crease B/H Placket (Single Fold) | Manual | 60 | 2 | 30.0 | 960 | Overstock of WIP |
| T22 | Top stitch B/H placket | Single Sewing MC | 40 | 1 | 40.0 | 720 | Less Stock of WIP |
| T23 | Sew Button Placket | Single Sewing MC | 25 | 1 | 25.0 | 1,152 | Overstock of WIP |
| T24 | Attach pocket | Manual | 20 | 1 | 20.0 | 1,440 | Overstock of WIP |
| T25 | Sewlabel at placket | Template Sewing MC | 50 | 2 | 25.0 | 1,152 | Overstock of WIP |
| T35 | Set front & back & mark neck for collar | Manual | 60 | 2 | 30.0 | 960 | Overstock of WIP |
| T36 | Shoulder attach | Single Sewing MC | 60 | 2 | 30.0 | 960 | Overstock of WIP |
| T37 | Shoulder top stitch | Single Sewing MC | 120 | 3 | 40.0 | 720 | Less Stock of WIP |
| T38 | Sleeve Attach | Single Sewing MC | 40 | 1 | 40.0 | 720 | Less Stock of WIP |
| T39 | Top stitch armhole | Single Sewing MC | 80 | 2 | 40.0 | 720 | Less Stock of WIP |
| T40 | Side Seam | Overlock MC | 110 | 3 | 36.7 | 785 | Less Stock of WIP |
| T41 | Collar Attach | Single Sewing MC | 30 | 1 | 30.0 | 960 | Overstock of WIP |
| T42 | Collar Close & Insert Label | Manual | 60 | 2 | 30.0 | 960 | Overstock of WIP |
| T43 | Cuff Attach & Close | Manual | 80 | 2 | 40.0 | 720 | Less Stock of WIP |
| T44 | Bottom Hem | Vertical Head Sewing M | 80 | 2 | 40.0 | 720 | Less Stock of WIP |
| T45 | Button Hold - Front Placket & Collar | Manual | 70 | 2 | 35.0 | 823 | Less Stock of WIP |
| T46 | Button Attach (Last Operation for finishing process) | Buttoning MC | 35 | 1 | 35.0 | 823 | Less Stock of WIP |
| | | | | | Through-put time: | 720 | (Result < Demand 900pcs) |

Table 5  What if all Online Production Cost Sheets for Style A

| Task | Online Qty (pcs) | Online Cost Rate (RMB/pc) | Offline Qty (pcs) | Offline Cost Rate (RMB/pc) | Total Cost (RMB) |
|---|---|---|---|---|---|
| T19 | 0 | 0.013 | 900 | 0.035 | 31.5 |
| T20 | 0 | 0.013 | 900 | 0.035 | 31.5 |
| T21 | 0 | 0.013 | 900 | 0.035 | 31.5 |
| T22 | 0 | 0.013 | 900 | 0.035 | 31.5 |
| T23 | 0 | 0.013 | 900 | 0.035 | 31.5 |
| T24 | 0 | 0.013 | 900 | 0.035 | 31.5 |
| T25 | 0 | 0.013 | 900 | 0.035 | 31.5 |
| T35 | 0 | 0.013 | 900 | 0.035 | 31.5 |
| T36 | 0 | 0.013 | 900 | 0.035 | 31.5 |
| T37 | 0 | 0.013 | 900 | 0.035 | 31.5 |
| T38 | 0 | 0.013 | 900 | 0.035 | 31.5 |
| T39 | 0 | 0.013 | 900 | 0.035 | 31.5 |
| T40 | 0 | 0.013 | 900 | 0.035 | 31.5 |
| T41 | 0 | 0.013 | 900 | 0.035 | 31.5 |
| T42 | 0 | 0.013 | 900 | 0.035 | 31.5 |
| T43 | 0 | 0.013 | 900 | 0.035 | 31.5 |
| T44 | 0 | 0.013 | 900 | 0.035 | 31.5 |
| T45 | 0 | 0.013 | 900 | 0.035 | 31.5 |
| T46 | 0 | 0.013 | 900 | 0.035 | 31.5 |
| Total Cost (RMB): | | | | | 598.5 |

Table 6  Optimized Production Cost Sheet for Style A

| Task | Online Qty (pcs) | Online Cost Rate (RMB/pc) | Offline Qty (pcs) | Offline Cost Rate (RMB/pc) | Total Cost (RMB) |
|---|---|---|---|---|---|
| T19 | 900 | 0.013 | 0 | 0.035 | 11.7 |
| T20 | 720 | 0.013 | 180 | 0.035 | 15.7 |
| T21 | 900 | 0.013 | 0 | 0.035 | 11.7 |
| T22 | 900 | 0.013 | 0 | 0.035 | 11.7 |
| T23 | 900 | 0.013 | 0 | 0.035 | 11.7 |
| T24 | 900 | 0.013 | 0 | 0.035 | 11.7 |
| T25 | 900 | 0.013 | 0 | 0.035 | 11.7 |
| T35 | 900 | 0.013 | 0 | 0.035 | 11.7 |
| T36 | 900 | 0.013 | 0 | 0.035 | 11.7 |
| T37 | 577 | 0.013 | 323 | 0.035 | 18.8 |
| T38 | 900 | 0.013 | 0 | 0.035 | 11.7 |
| T39 | 900 | 0.013 | 0 | 0.035 | 11.7 |
| T40 | 785 | 0.013 | 115 | 0.035 | 14.2 |
| T41 | 900 | 0.013 | 0 | 0.035 | 11.7 |
| T42 | 900 | 0.013 | 0 | 0.035 | 11.7 |
| T43 | 900 | 0.013 | 0 | 0.035 | 11.7 |
| T44 | 720 | 0.013 | 180 | 0.035 | 15.7 |
| T45 | 900 | 0.013 | 0 | 0.035 | 11.7 |
| T46 | 822 | 0.013 | 78 | 0.035 | 13.4 |
| Total Cost (RMB): | | | | | 241.6 |

Overall, the findings confirm that implementing MILP for line balancing not only enhances the alignment of production capabilities with market demand but also significantly reduces labour costs. This optimization strategy provides a compelling case for manufacturers in the garment industry to adopt advanced mathematical modelling techniques to improve operational efficiency and achieve substantial cost savings. The results validate the expected outcomes of MILP, reinforcing its value as a strategic tool for optimizing production processes in the modern garment industry.

B. *Comparative Analysis*
- Throughput Enhancement:

  The optimization plan addressed the throughput inefficiency, aligning production output more closely with the daily demand of 900 pieces.

- Workload Balancing:

  By systematically formulating the problem with appropriate variables, objective functions, and constraints, the MILP approach balanced workloads across both online and offline production environments.

- Financial Impact:

  The substantial reduction in daily costs underscores the effectiveness of the MILP approach in optimizing resource allocation and minimizing waste.

The integration of Mixed-Integer Linear Programming (MILP) with Enterprise Resource Planning (ERP) systems for line balancing represents a forward-thinking approach to optimizing production processes. By developing proprietary software tailored to this integration, organizations can achieve several strategic advantages:

C. *Enhanced Operational Efficiency*
Seamless Integration:

Combining MILP with ERP systems allows for real-time data exchange and decision-making, enabling more responsive and adaptive production line management.

Automated Optimization:

The software can automate the optimization of task assignments and resource allocation, reducing manual intervention and minimizing human error.

D. *Improved Decision-Making*
- Data-Driven Insights:

  Leveraging the comprehensive data capabilities of ERP systems, the integrated software can provide actionable insights and predictive analytics to support strategic planning and operational adjustments.

- Scenario Analysis:

  The ability to simulate various production scenarios and their outcomes empowers managers to make informed decisions that align with business objectives and market demands.

E. *Cost and Resource Management*
- Cost Reduction:

  By optimizing line balancing through MILP, the software can help reduce both direct and indirect labour costs, as demonstrated by the significant cost savings achieved in previous analyses.

- Resource Utilization:

  Enhanced visibility into resource availability and utilization ensures that production processes are aligned with demand, minimizing waste and improving overall efficiency.

F. *Competitive Advantage*
- Customization and Flexibility:

  Developing proprietary software allows for customization to meet specific organizational needs and industry requirements, providing a competitive edge over standardized solutions.

- Scalability:

  The software can be designed to scale with the organization, accommodating growth and evolving production complexities without compromising performance.

### G. Considerations for Development

- Technical Expertise:

  Developing such software requires expertise in both MILP and ERP systems, as well as a deep understanding of the production processes and industry-specific challenges.

- Investment and Resources:

  Significant investment in terms of time, financial resources, and personnel will be necessary to develop, implement, and maintain the software.

- Change Management:

  Successful implementation will require effective change management strategies to ensure user adoption and integration into existing workflows.

In conclusion, the development of new software to integrate MILP with ERP systems for line balancing holds significant promise for enhancing production efficiency and strategic decision-making. By investing in this initiative, organizations can position themselves at the forefront of innovation in manufacturing processes, driving long-term success and competitiveness in the market.

This article has shown the successful case for applying the new MILP technology for Line Balancing with ERP in the Modern Garment Industry, referring to the work of Prof Ray WM Kong [9].

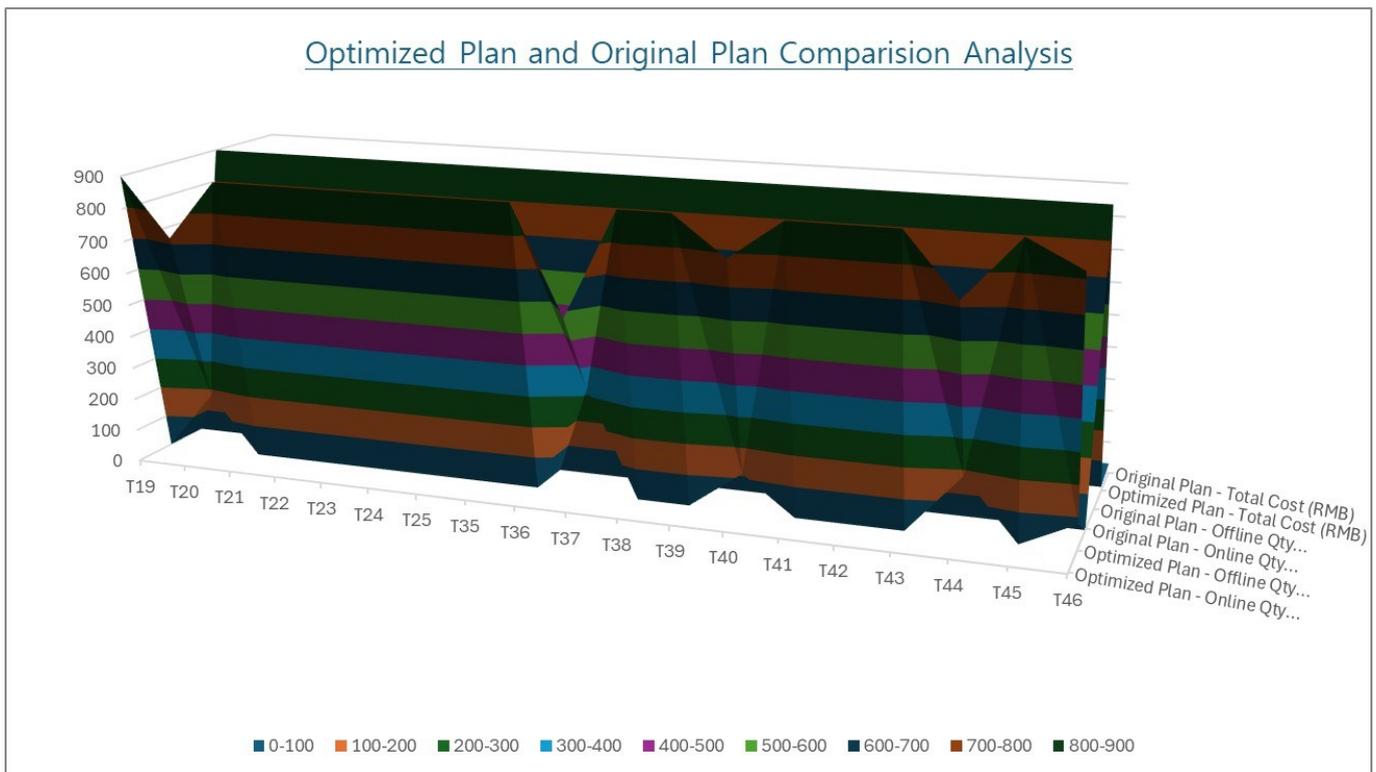

Fig 9 Optimized Plan and Original Plan Comparison Analysis Sheet

**Reference**

[1]. Ray Wai Man Kong, Theodore Ho Tin Kong, Tianxu Huang (2024) Lean Methodology for Garment Modernization, IJEDR – International Journal of Engineering Development and Research (www.IJEDR.org), ISSN:2321-9939, Vol.12, Issue 4, page 14-29, October-2024,
(https://rjwave.org/IJEDR/papers/IJEDR2404002.pdf)

[2]. Bongomin O, Mwasiagi JI, Nganyi EO, Nibikora I. (2020) A complex garment assembly line balancing using simulation-based optimization. Engineering Reports. 2020; 2:e12258
 (https://doi.org/10.1002/eng2.12258)

[3]. Gary Yu-Hsin Chen, Ping-Shun Chen, Jr-Fong Dang, Sung-Lien Kang, Li-Jen Cheng (2021), Applying Meta-Heuristics Algorithm to Solve Assembly Line, International Journal of Computational Intelligence Systems, Vol. 14(1), 2021, pp. 1438–1450
   (https://doi.org/10.2991/ijcis.d.210420.002)

[4]. R. Chen, C. Liang, D. Gu, J.Y.-T. Leung, A multi-objective model for multi-project scheduling and multi-skilled staff assignment for IT product development considering competency evolution, Int. J. Prod. Res. 55 (2017), 6207–6234.

[5]. Prabir Jana, Different Machine Layouts in Sewing Line, Apparel Resouce, web article in the Apparel Resource (2009),(https://apparelresources.com/business-news/manufacturing/basics-of-machine-layout-in-sewing-line)

[6]. Kong, R.W., Liu, M., & Kong, T.H. (2024), Design and Experimental Study of Vacuum Suction Grabbing Technology to Grasp Fabric Piece, Open Access Library

Journal, v11, pg 1-17, (DOI: 10.4236/oalib.1112292), ArXiv, abs/2408.09504. (https://www.oalib.com/research/6838637)

[7]. Ray Wai Man Kong, Theodore Ho Tin Kong, Miao Yi (2024) Design a New Pulling Gear for the Automated Pant Bottom Hem Sewing Machine, International Research Journal of Modernization in Engineering Technology and Science, Vol. 06, No. 11, 11.2024, p. 3067-3077. DOI: 10.56726/IRJMETS64156 (https://www.doi.org/10.56726/IRJMETS64156)

[8]. Kong, R.W., Kong, T.H., & Huang, T. (2024). Lean Methodology for Garment Modernization. ArXiv, abs/2410.07705. (https://arxiv.org/abs/2410.07705)

[9]. Prof. Dr. Ray Wai Man Kong, Ding Ning, & Theodore Ho Tin Kong. (2025). Line Balancing in the Modern Garment Industry. International Journal of Mechanical and Industrial Technology, 12(2), 60–72. (https://doi.org/10.5281/zenodo.14800724)

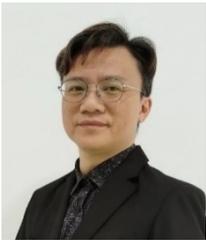

**Ray Wai Man Kong** (Senior Member of IEEE, member of IET, member of Public Administration Association) Hong Kong, China. He received a Bachelor of General Study degree from the Open University of Hong Kong, Hong Kong in 1995. He received an MSc degree in Automation Systems and Engineering and an Engineering Doctorate from the City University of Hong Kong, Hong Kong in 1998 and 2008 respectively.

From 2005 to 2013, he was the operations director with Automated Manufacturing Limited, Hong Kong. From 2020 to 2021, he was the Chief Operating Officer (COO) of Wah Ming Optical Manufactory Ltd, Hong Kong. He is a modernization director with Eagle Nice (International) Holdings Limited, Hong Kong. He is an Adjunct Professor of the System Engineering Department at the City University of Hong Kong, Hong Kong. He published more articles in international journals on robotic technology, grippers, automation and intelligent manufacturing. His research interests focus on intelligent manufacturing, automation, maglev technology, robotics, mechanical engineering, electronics, and system engineering for industrial factories.

Prof. Dr. Kong Wai Man, Ray is Vice President of CityU Engineering Doctorate Society, Hong Kong and a chairman of the Intelligent Manufacturing Committee of the Doctors Think Tank Academy, Hong Kong. He has published more intellectual properties and patents in China.

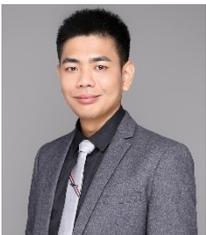

**James Ding Ning** （IEEE member, Engineering Doctorate Student）Hong Kong, China. He received a Bachelor of Electronic and Information Engineering degree from Shenzhen University, Guangdong Province, China in 2006. He received an MSc degree in Electronic Engineering from the University of Sheffield, United Kingdom in 2007 with upper second honours. He is a part-time Engineering Doctorate Student of the System Engineering Department at the City University of Hong Kong.

Mr. Ning has more than 16 years of industry experience focused on semiconductor and intelligence instrument design. He is the first author of more than 10 patents. He is currently working at Huawei Investment Limited Corporation, Hong Kong as a Senior Engineer. His research interest is focused on computer architecture design, CPU Design, integrated circuit design, communication standard protocol and industry management. He has extensive work experience in both Shenzhen and Hong Kong, with his professional footprint spanning across the Greater Bay Area.

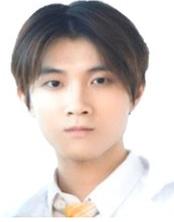

**Theodore Ho Tin Kong** (MIEAust, Engineers Australia and MIEEE) received his Bachelor of Engineering (Honours) in mechanical and aerospace engineering from The University of Adelaide, Australia, in 2018. He then earned a Master of Science in aeronautical engineering (mechanical) from HKUST - Hong Kong University of Science and Technology, Hong Kong, in 2019.

He began his career as a Thermal (Mechanical) Engineer at ASM Pacific Technology Limited in Hong Kong, where he worked from 2019 to 2022. Currently, he is a Thermal-Acoustic (Mechanical) Design Engineer at Intel Corporation in Toronto, Canada. His research interests include mechanical design, thermal management and heat transfer, and acoustic and flow performance optimization. He is proficient in FEA, CFD, thermal simulation, and analysis, and has experience in designing machines from module to heavy mechanical level design.